\documentclass[12pt,leqno]{article}
\def\s{\dot{s}} 
\def\w{\dot{w}}
\def\v{{\rm v}}
 
\def\int{\mathbb{Z}}

\def\Ue{{\cal U}_{\varepsilon}({\mathfrak g})}

\def\O{{\cal O}}

\def\proof{{\bf Proof. }}

\def\pf{\proof}

\def\V{{\cal V}}
\usepackage{times}
\usepackage{graphics}
\usepackage{amssymb}
\usepackage{amscd}
\usepackage{amsmath}
\input xy 
\xyoption{all}
\title{A classification of spherical conjugacy classes in good characteristic}
\newtheorem{theorem}{Theorem}[section]
\newtheorem{lemma}[theorem]{Lemma}

\newtheorem{definition}[theorem]{Definition}
\newtheorem{remark}[theorem]{Remark}

\author{Giovanna Carnovale\\
Dipartimento di Matematica Pura ed Applicata\\
Torre Archimede - via Trieste 63 - 35121 Padova - Italy\\
email: carnoval@math.unipd.it }
\date{}
\begin{document}
\maketitle
\begin{abstract}We classify spherical conjugacy classes in a simple
  algebraic group over an algebraically closed field of good, odd characteristic.
\end{abstract}



\section*{Introduction}

When studying a transitive action of a group $G$, it is particularly
interesting to understand when a given subgroup $B$ of $G$ acts with
finitely many orbits. A particular case of this situation in  the
theory of algebraic groups is when $B$ is a Borel
subgroup of a connected algebraic group $G$. The $G$-spaces for which
$B$ acts with finitely many orbits are the so-called spherical
homogeneous spaces and they include important examples such as the
flag variety $G/B$ and symmetric spaces. They are precisely those
$G$-spaces for which the $B$-action has a dense orbit for the Zariski
topology. One may want to understand when homogeneous spaces which are
relevant in algebraic Lie theory, such as nilpotent orbits in ${\rm Lie}(G)$ and conjugacy classes in $G$, are spherical.
Spherical nilpotent orbits have been classified in \cite{pany,pany2}
when the base field is ${\mathbb C}$ and in \cite{FR} when it is an
algebraically closed field of good characteristic:
they are precisely the orbits of type $rA_1$ for $r\geq0$ in the simply-laced case and of type $rA_1+s\tilde{A}_1$ for $r,s\geq0$ in the multiply-laced case.
As for conjugacy classes, it is natural to use the interplay with the Bruhat decomposition, which has proven to be a fruitful tool in the past. For instance, it is essential in the description of regular classes (\cite{cross-section}) whose intersection with Bruhat cells is still subject of current research (\cite{EG,EG2}). This approach has brought 
to two characterizations of the spherical conjugacy
classes in a connected, reductive algebraic group $G$ over an
algebraically closed field of zero or good, odd characterisitc
(\cite{ccc,gio,gio2}). The first one is given through a
formula relating the dimension of a class and the Weyl
group element corresponding to a suitable Bruhat cell intersecting the class. The second one states that spherical conjugacy classes are
exactly those classes intersecting only Bruhat cells corresponding to
involutions in the Weyl group of $G$. These characterizations can be
used in order to give a complete and exhaustive list of the spherical classes of $G$. In a simple algebraic group over ${\mathbb C}$ they have been classified in \cite{ccc}, making use of the classification of spherical nilpotent orbits. Spherical classes in type $G_2$ in good characteristic have been classified in  \cite{gio2}. 

The present paper completes the picture, as it classifies spherical
classes in good, odd characteristic. Contrarily to \cite{ccc}, the present work is independent of the
classification of spherical unipotent conjugacy classes existing in
the literature. Since Springer corrispondence holds in good characteristic,
it provides an elementary classification of
spherical nilpotent orbits alternative to \cite{FR}, where
Kempf-Russeau theory is involved and where the aid of a computer programme is needed to deal with the exceptional types.
The crucial tools in our method are just those conditions in the
characterizations in \cite{ccc,gio,gio2} whose proofs are general and
rather short.  
The arguments used for this classification can also be
transferred to the characteristic zero situation, providing an
alternative, elementary approach to \cite{pany,pany2}, although by case-by-case considerations. 

After fixing notations and recalling basic notions in \S
\ref{preliminaries},  we introduce spherical conjugacy classes and
their characterizations in \S \ref{chara}. We proceed providing
the list of spherical conjugacy classes through a case-by-case
analysis in \S \ref{classi}.

The result is as when the base field is ${\mathbb C}$:
In the simply-laced case spherical conjugacy classes are either
semisimple or unipotent and the semisimple ones are
all symmetric spaces if $G$ is simply-connected.
In type $G_2$ spherical conjugacy classes are again either semisimple or unipotent but, as in types $B_n$ and $C_n$, there are spherical semisimple classes that are not symmetric. Just as in other situations involving spherical homogeneous spaces
(for example, in the description of maximal spherical ideals of Borel
subalgebras \cite{PR}) the doubly-laced case is
slightly more involved. The new phenomenon in the present situation is
that there appear spherical classes
that are neither semisimple nor unipotent. 

\section{Notation}\label{preliminaries}

Unless otherwise stated $G$ will denote
a simple algebraic group over an algebraically closed field
$k$ of good odd characteristic (\cite[\S I.4]{131}). When we consider an
integer as an element in $k$ we shall mean its image
in the prime field of $k$. We shall denote by  $\Phi$ the root system 
relative to a fixed Borel subgroup $B$ and a maximal torus $T$ of $G$; by
$\Delta=\{\alpha_1, \dots, \alpha_n\}$
the corresponding set of simple roots and by $\Phi^+$ the set of positive roots. We shall use the numbering of
the simple roots in \cite[Planches I-IX]{bourbaki}. The highest positive root will be denoted by $\beta_1$. 
For generators of $T$ we shall adopt the notation in \cite[Lemma 19]{lectures}. We will put $W=N(T)/T$ and
$s_{\alpha}$ will indicate the reflection corresponding to the root $\alpha$. Given an element $w\in W$ we shall denote by $\w$
a representative of $w$ in $N(T)$. The maximal unipotent subgroup of $B$ will be denoted by $U$.
For a root $\alpha$ the elements of the associated root subgroup
$X_\alpha$ will be denoted by $x_{\alpha}(t)$.
%
%
We shall frequently use that $x_{-\alpha}(t)\in B s_\alpha B$ and that
if $\alpha\in\Delta$ and $w\alpha\in\Phi^+$ then $BwBs_\alpha
B=Bws_\alpha B$ whereas if $w\alpha\in-\Phi^+$ then $BwBs_\alpha
B\subset BwB\cup Bws_\alpha B$ (\cite[Lemma 8.3.7]{springer}).

Given an element $x\in G$ we shall denote by $\O_x$ the conjugacy
class of $x$ in $G$ and by $H_x$ the centralizer of $x$ in $H\leq
G$.  For the dimension of unipotent conjugacy classes
in arbitrary good characteristic we refer to \cite[Chapter
  13]{carter-finite} and \cite[Theorem 2.6]{premet}. 

\section{Characterizations through the Bruhat decomposition}\label{chara}

In this section we shall introduce the characterizations of spherical conjugacy classes that we will use in the sequel. 

\begin{definition}\label{sferica} Let $G$ be a connected reductive algebraic
  group. A homogeneous space $G/H$ is {\em spherical} if it has a dense orbit for some Borel subgroup of $G$.
\end{definition}

By an abuse of notation we shall say that $g\in G$ is spherical if its class
$\O_g\simeq G/G_g$ is so.\\
For a conjugacy class $\O$ in $G$, we shall denote by $\V$ the set of
$B$-orbits in $\O$. It is well-known (\cite{Bri,Vin} in characteristic $0$,
\cite{gross,knop} in positive characteristic)
that $\O$ is a spherical conjugacy class if and only if $\V$ is finite. 

Since $G=\bigcup_{w\in W}BwB$, for every class $\O$  there is a natural map
$\phi\colon \V\to W$ associating to
$v\in\V$ the element $w$ in the Weyl group of $G$ for which
$v\subset BwB$. Besides, there is a unique $w\in W$ for which $BwB
\cap\O$ is dense in $\O$ and this element is maximal in ${\rm
  Im}(\phi)$ with respect to the Bruhat ordering (\cite[p. 32]{ccc}). We shall denote such an element by $w_\O$.

\smallskip

There are two characterizations of spherical classes in
$G$. Let $\ell$ denote the usual length function on $W$ and let ${\rm rk}(1-w)$ denote the rank of the operator $1-w$ in the geometric representation of $W$.

\begin{theorem}\label{dime}(\cite[Theorem 25]{ccc}, \cite[Theorem
    4.4]{gio}) A class $\O$ in a connected reductive algebraic group
      $G$ over an algebraically closed field of zero or good odd characteristic is spherical if and only if there exists $v$ in $\V$ such that $\ell(\phi(v))+{\rm rk}(1-\phi(v))=\dim\O$. If this is the case, $v$ is the dense $B$-orbit and $\phi(v)=w_\O$.
\end{theorem} 

\begin{theorem}\label{invo}(\cite[Theorem 2.7]{gio},\cite[Theorem
    5.7]{gio2}) A class $\O$ in  a connected reductive algebraic group
  $G$ over an algebraically closed field of zero or odd, good
  characteristic is spherical if and only if ${\rm Im}(\phi)$ contains only involutions in $W$. 
\end{theorem}

\begin{remark}\label{regular}{\rm Regular classes
 in a reductive algebraic group whose semisimple quotient is not of
 type $rA_1$ cannot be
 spherical. Indeed, by \cite[Theorem 8.1]{cross-section}, regular
 classes intersect  Bruhat cells corresponding to Coxeter elements.} 
\end{remark}

\begin{remark}\label{same}{\rm Let $g$ and $x$ be elements in $G$ with $G_x=G_g$. Then $\O_g\simeq G/G_g$ is spherical if and only if $\O_x\simeq G/G_x$ is so.
In particular, if $g^2\in Z(G)$ then $\O_x$ is a symmetric space, hence spherical by \cite[Corollary 4.3]{results}.}
\end{remark}

Let $g\in B$ with Jordan decomposition $g=su\in TU$. Then $u\in
G_s^\circ$ and $G_g=G_s\cap G_u$. Therefore, if  $BG_g$ is dense in $G$
then $BG_s$ and $BG_u$ are dense in $G$. In other words, if $\O_g$ is spherical then $\O_s$ and $\O_u$ are also spherical. We can refine this argument.

\begin{lemma}\label{mixed}Let $g\in B$ with Jordan decomposition
  $g=su\in TU$. If $\O_g$ is spherical then $\O_s$ and $\O_u$ are spherical in $G$ and the class $\O'_u$ of $u$ in $G_s^\circ$ is spherical in $G_s^\circ$. 
\end{lemma}
\pf The subgroup $B_1=B\cap G_s^\circ$ is a Borel subgroup for
$G_s^\circ$ containing $T$. If $x=b_1 \w b_2\in B_1 (N(T)\cap G_s^\circ) B_1$ is a Bruhat
decomposition of $x\in G_s^\circ$ then it is also a Bruhat
decomposition of $x$ in $G$. If $\v\in\O'_u$
with $\v=b_1\w b_2$ then $s\v=sb_1 \w b_2\in B N(T) B\cap \O_g$ forcing $\w^2\in T$ by Theorem \ref{invo} applied to $G$. By Theorem \ref{invo} applied to $G_s^\circ$ we deduce that $\O'_u$ is spherical therein.\hfill$\Box$

\section{The classification}\label{classi}

We aim at a classification of spherical conjugacy classes in good odd
characteristic since the classification for $k={\mathbb
  C}$ in \cite{ccc} holds for every
algebraically closed field of characteristic zero. The property
of being spherical for $G/H$ depends only on ${\rm Lie}(G)$ and ${\rm Lie}(H)$
so the classification will not distinguish groups  isogenous to $G$ and it will be up
to a central element in $G$. The main tools will be the
sufficient condition in Theorem \ref{dime} (\cite[Theorem 5]{ccc}) and the
necessary condition in Theorem \ref{invo} (\cite[Theorem 2.7]{gio}). If $G$ is of type $G_2$ the classification in good
characteristic is given in \cite[\S 2.1]{gio2} and we provide it here
for the sake of completeness. 

\subsection{Type $G_2$}\label{g2}

\begin{theorem}Let $G$ be of type $G_2$. The spherical
  classes are either semisimple or unipotent. The semisimple ones
  are 
represented by $h_{\alpha_1}(-1)$ and
  $h_{\alpha_1}(\zeta)$ for $\zeta$ a fixed primitive third root of
  $1$. The unipotent ones are those of type $A_1$ and $\tilde{A}_1$.
\end{theorem}

\subsection{Type $A_n$}

\begin{theorem}Let $G=SL_{n+1}(k)$. If $n=1$ all
  classes are spherical. If $n\geq2$ the spherical classes
  are either semisimple or unipotent. The semisimple ones are those corresponding to matrices with at
  most two distinct eigenvalues and they are all symmetric
  spaces. The unipotent ones are those associated with partitions of
  type $(2^m,1^{n+1-2m})$ for $m=1,\ldots,[\frac{n+1}{2}]$.
\end{theorem}
\pf If $n=1$ all Bruhat cells correspond to involutions in $W$ so
every class is spherical by Theorem \ref{invo}.  

Let $G$ be $SL_{n+1}(k)$ with $n\geq2$, let $B, T, U$ be as usual and let $\O=\O_u$ be a unipotent class. By Jordan theory we may assume that $u=x_{-\alpha_1}(c_1)\cdots x_{-\alpha_n}(c_n)$ with $c_i\in\{0,\,1\}$.  By \cite[Lemma 8.1.4(i), Lemma 8.3.7]{springer} this element lies in $Bs_1^{c_1}\cdots s_n^{c_n}B$. This cell corresponds to an involution only if $c_ic_{i+1}=0$ for all $i=1,\ldots, n-1$. Theorem \ref{invo} implies that if $\O_u$ is spherical its associated partition is of type $(2^m,1^{n+1-2m})$. Conversely, let $\O$ be a unipotent class corresponding to $(2^j,1^{n+1-2j})$, with $2j\leq n+1$. Let $\beta_i=\alpha_i+\cdots+\alpha_{n-i+1}$ for $i=1,\ldots,j$.
The element $x_{-\beta_1}(1)\cdots x_{-\beta_j}(1)$ lies in $\O$ and 
 its $B$-orbit satisfies the condition in Theorem \ref{dime} so $\O$ is spherical. 

Let $s={\rm diag}(\lambda_1 I_{n_1},\lambda_2
 I_{n_2},\ldots,\,\lambda_lI_{n_l})\in T$. If $l$ is greater than $2$ then $s$ is conjugate to some matrix
$t={\rm diag}(\lambda_1,\lambda_2,\lambda_3, t_1)$ for some invertible
 diagonal submatrix $t_1$. Then $t$ lies in the reductive subgroup  $H=\langle
 T,X_{\pm\alpha_1},X_{\pm\alpha_2}\rangle$ and it is regular
 therein. Besides, $B_1=B\cap H$ is a Borel subgroup of $H$ containing $T$. 
By Remark \ref{regular} some conjugate of $t$ 
lies in a Bruhat cell of $H$ that does not correspond to an
involution. Therefore, if $\O_s$ is spherical semisimple, $s$ has at most $2$ eigenvalues.
%
Conversely,
 suppose that $s\in T$ has $2$ eigenvalues. We may assume
 $s={\rm diag}(\lambda I_{m},\mu I_{n+1-m})$. Let $\zeta$ be a
 primitive $2(n+1)$-th root of unity if $n+1-m$ is odd and let
 $\zeta=1$ if $n+1-m$ is even. Let also $s_0={\rm diag}(\zeta
 I_{m},-\zeta I_{n+1-m})$. Then $s_0^2\in Z(G)$ and  $G_s=G_{s_0}$. By
 Remark \ref{same} the class $\O_s$ is symmetric.

We will now show that there is no spherical element with Jordan
decomposition  $x=su$ with $s\not\in Z(G)$ and $u\neq1$. Were this the
case, we could assume that $s={\rm diag}(\lambda I_m,\mu I_{n+1-m})$
with $m\geq2$ and that $x$ is conjugate to some $y=\begin{pmatrix}
\lambda u_1&\\
&\lambda&\\
&1&\lambda&\\
&&&\mu u_2
\end{pmatrix}$
with $u_1\in X_{\alpha_1}\cdots X_{\alpha_{m-3}}$ and $u_2\in
X_{\alpha_{m+1}}\cdots X_{\alpha_n}$. Then
$y'=\begin{pmatrix}
\lambda u_1\\
&\lambda\\
&1&\lambda\\
&&1&\mu u_2\\
\end{pmatrix}\in \O_x\cap Bs_{m-1}s_m B$. 
By Theorem \ref{invo} we have the statement.\hfill$\Box$

\subsection{Type $C_n$}
\begin{theorem}\label{cn}Let $G=Sp_{2n}(k)$ for $n\geq2$. The spherical semisimple
  classes are represented up to a sign by $\sigma_l={\rm
  diag}(-I_l,I_{n-l},-I_l, I_{n-l})$ for $l=1,\ldots,\left[\frac{n}{2}\right]$; $a_\lambda={\rm diag}(\lambda I_n,
  \lambda^{-1} I_n)$ and $c_\lambda={\rm
  diag}(\lambda, I_{n-1},\lambda^{-1}, I_{n-1})$ for $\lambda^2\neq0,1$. The unipotent ones are those whose associated partition is of type $(2^m, 1^{2n-2m})$ for $m=1,\ldots,n$. The spherical
  classes that are neither semisimple nor unipotent are
  represented up to a sign by the elements $\sigma_lu$
  where $u\in G_{\sigma_l}\cong Sp_{2l}(k)\times Sp_{2n-2l}(k)$ is unipotent and corresponds to the partition $(2, 1^{2n-2})$.  
\end{theorem}
\pf We view $G$ as the subgroup of $GL_{2n}(k)$ of matrices preserving the bilinear form associated with the matrix $\begin{pmatrix}0&I\cr
-I&0\cr\end{pmatrix}$ with respect to the canonical basis of $k^{2n}$. 
We choose $B$ as the subgroup of $G$ of matrices of the form $\begin{pmatrix}A&AX\cr
0&^t\!A^{-1}\cr\end{pmatrix}$ where $A$ is an invertible upper
triangular matrix, $^t\!A^{-1}$ is its inverse transpose and $X$ is a
symmetric matrix. The torus $T$ is the subgroup of diagonal matrices in $B$.

Let $s\in T$. If $s$ has at least $4$ eigenvalues then it
is either conjugate to a matrix of type $s'={\rm
  diag}(r,\lambda,\mu,r^{-1},\lambda^{-1},\mu^{-1})$ for some $\lambda,\,\mu\in k$ with
$\lambda\not\in\{\lambda^{-1},\mu,\mu^{-1}\}$ and $\mu^2\neq1$ and
some invertible diagonal submatrix $r$, or to a matrix of type $s''={\rm
  diag}(\lambda,1,-1,t,\lambda^{-1},1,-1,t^{-1})$ for some
invertible diagonal submatrix $t$ and some $\lambda\in k$ with $\lambda^2\neq1$. 
The element $s'$ is regular in $H=\langle T,\,X_{\pm\alpha_{n-1}},\,X_{\pm \alpha_n}\rangle$ whereas  $s''$ is regular in $H'=\langle T,\,X_{\pm\alpha_1},\,
X_{\pm\alpha_2}\rangle$. By Remark \ref{regular} applied to $H$ and $H'$ 
 the class $\O_s$ cannot be spherical. 

Let now $s$ have exactly $3$ eigenvalues. Then one of these is $\pm1$
and it is not restrictive to assume that they are
$\lambda,\lambda^{-1},1$ with $\lambda^2\neq1$. If the
multiplicity of $\lambda^{\pm1}$ is greater than $1$, $s$ is conjugate to some $r'={\rm
  diag}(\lambda,\lambda^{-1},1,r_1,\lambda^{-1},\lambda, 1, r_1^{-1})$
with $r_1$ an invertible, diagonal submatrix. The element $r'$ lies and is regular in the subgroup $H'$ described above. By Remark \ref{regular} the class $\O_s$ cannot be spherical.  
%
On the other hand if  $s$ has exactly $3$ eigenvalues
$\lambda^{\pm1}, 1$ (up to a sign) with the multiplicity of
$\lambda^{\pm1}$ equal to $1$, then $\O_s$ is spherical. Indeed, the
representative of its class in \cite[Theorem 15]{ccc} works also in odd characteristic and its $B$-orbit satisfies the condition of Theorem \ref{dime}.
%

Let $s\in T$ have exactly $2$ eigenvalues. They are either $1,\,-1$ so that $\O_s$ is symmetric, or $\lambda, \lambda^{-1}$ for $\lambda^2\neq 1$. In the latter case $s$ is conjugate to ${\rm diag}(\lambda I_n, \lambda^{-1}I_n)$ and its centralizer coincides with the centralizer of $s_0={\rm diag}(\zeta I_n, -\zeta I_n)$ for $\zeta$ a primitive fourth root of $1$. Since $\O_{s_0}$ is a symmetric space, $\O_s$ is spherical by Remark \ref{same}. 

Unipotent classes in $G$ are parametrized through Jordan theory by partitions where odd terms occur pairwise. 
Let $\O_u$ be a unipotent class and let $\underline{\lambda}$ be its associated partition.
Let $\underline{\mu}=(\mu_1,\ldots,\mu_l)$ be obtained by taking a representative of each term occurring pairwise in $\underline{\lambda}$ and let $\underline{\nu}=(\nu_1,\ldots,\nu_m)$ be obtained by taking the (even) terms without repetition in $\underline{\lambda}$, so that $2n=|\underline{\nu}|+2|\underline{\mu}|$. A representative $u'$ of $\O_u$ can be taken in the subgroup isomorphic to $Sp_{2\mu_1}(k)\times\cdots\times Sp_{2\mu_l}(k)\times Sp_{\nu_1}(k)\times\cdots \times Sp_{\nu_m}(k)$ obtained by repeating the immersion of  $Sp_{2d_1}(k)\times Sp_{2d_2}(k)$ into $Sp_{2(d_1+d_2)}(k)$ given by $\left(\begin{pmatrix}A_1&B_1\\C_1& D_1\end{pmatrix}, \begin{pmatrix}A_2& B_2\\
C_2& D_2\end{pmatrix}\right)\mapsto\begin{pmatrix}
A_1&&B_1\\
&A_2&&B_2\\
C_1&&D_1\\
&C_2&&D_2
\end{pmatrix}$. The component of $u'$ in $Sp_{\nu_j}(k)$ corresponds to the partition $(\nu_j)$. Besides, the component of $u'$ in $Sp_{2\mu_i}(k)$ can be taken to lie and be regular in the subgroup isomorphic to $SL_{\mu_i}(k)$ obtained through the immersion $M\mapsto {\rm diag}(M, ^t\!M^{-1})$. In particular, the choices of Borel subgroups we made are all compatible and $u'$ is regular in a subgroup of type $A_{\mu_1-1}\times\cdots\times A_{\mu_l-1}\times C_{\nu_1}\times\cdots\times C_{\nu_m}$. Therefore, if $u$ is spherical, we necessarily have $\mu_i=2$ and  $\nu_j=1$ for every $i$ and $j$. 
Conversely, let $\underline{\lambda}=(2^j,1^{2n-2j})$ and let $\O_j$ be the unipotent class associated with $\underline{\lambda}$. Then for $\beta_q=2\alpha_q+\cdots+2\alpha_{n-1}+\alpha_n$ for $q=1,\ldots, n$ the $B$-orbit of the element $x_{-\beta_1}(1)\cdots x_{-\beta_j}(1)$ satisfies the condition in Theorem \ref{dime} for $\O_j$ (cfr. \cite[Theorem 12]{ccc}).

Let $g=su$ be the Jordan decomposition of a spherical element in $G$ with $s\not\in Z(G)$ and $u\neq 1$. Then $\O_s$ is spherical so we may assume $s=a_\lambda, c_\lambda$ or $\sigma_l$. The case $s=a_\lambda$ is ruled out because we would have $\dim O_{a_\lambda u}>\dim \O_{a_\lambda}=\dim B$. Let us assume that $s=c_\lambda$. 
Then $u\in G_s\cong k^*\times Sp_{2n-2}(k)$ and it is spherical therein so
 it corresponds to a partition of the form $(2^m,1^{2n-2-2m})$. 
The class $\O_{c_\lambda u}$ may be represented by a matrix of the form $\begin{pmatrix}A\\
^t\!A^{-1}X&^t\!A^{-1}\end{pmatrix}$ where  
$A=\begin{pmatrix}\lambda\\
1&1\\
&&I_{n-2}\\
\end{pmatrix}$ and $X={\rm diag}(0,I_m,0_{n-m-1})$.
Such a representative is contained in $T X_{-\alpha_1}X_{-\beta_2}\cdots X_{-\beta_{m+1}}$ with notation as above, and it  
lies in $Bs_1s_{\beta_2}\cdots s_{\beta_{m+1}}B$. By Theorem \ref{invo}  this is not possible. It follows that we necessarily have $s=\sigma_l$ for some $l$. 

In this case $G_s$ is a subgroup isomorphic to $Sp_{2l}(k)\times
Sp_{2n-2l}(k)$ and $u=(u_1,u_2)$ lies in $G_s$ and it is spherical
therein. Then  $u_1$ and $u_2$ are spherical in the respective
components. We claim that $u_1$ and $u_2$ correspond to partitions
with no repeated terms. If
$\underline{\lambda}=(2,2,\underline{\lambda}')$ were the partition
associated with $u_2$ then $\O_u$ would contain a matrix of the form
$\begin{pmatrix}
A_1&&&A_1X_1\\
\begin{array}{cc}
&1\\
\\
\end{array}&\begin{array}{cc}
1\\
1&1
\end{array}\\
&&A_2&&&A_2X_2\\
&&&^t\!A_1^{-1}&\begin{array}{cc}
&\\
1\phantom{11}&\\
\end{array}\\
&&&&\begin{array}{cc}
1&-1\\
&1
\end{array}\\
&&&&&^t\!A_2^{-1}
\end{pmatrix}
$
where $\begin{pmatrix}
A_1&A_1X_1\\
&^t\!A_1^{-1}\end{pmatrix}$ represents $-u_1$, $\begin{pmatrix}
A_2& A_2X_2\\
&^t\!A_2^{-1}\end{pmatrix}$  represents a unipotent element
corresponding to $\underline{\lambda'}$ and $A_1,A_2$ are upper
triangular. Such a representative lies in $B s_ls_{l+1}B$ leading to a
contradiction. The case in which the partition associated with $u_1$
has repeated terms can be treated similarly. We claim that $u$ can have only one nontrivial component. If this were
not the case, the class would be represented by a matrix of the form 
$\begin{pmatrix} 
A\\
^t\!A^{-1}X&^t\!A^{-1}
\end{pmatrix}$ where $A=\begin{pmatrix}
-I_{l-1}\\
&-1\\
&1&1\\
&&&I_{n-l-1}
\end{pmatrix}$ and $X={\rm diag}(0_{l},1,0_{n-l-1})$. Such a representative lies in $TX_{-\alpha_l}X_{-\beta_{l+1}}$ and its Bruhat cell corresponds to $s_ls_{\beta_{l+1}}$ which is not an involution.  
%
Conversely, for all classes of type $\sigma_lu$ with $u\in G_{\sigma_l}$ corresponding to the partition $(2, 1^{2n-2})$ the representative in \cite[Theorem 21]{ccc} is defined in odd characteristic and its $B$-orbit satisfies the condition of Theorem \ref{dime}. \hfill$\Box$

\subsection{Type $D_n$}
\begin{theorem}\label{dn}Let $G=SO_{2n}(k)$ for $n\ge4$. The spherical  classes in $G$ are either semisimple or unipotent. The semisimple
ones are represented, up to a central element and up to the automorphism arising from an automorphism of the Dynkin diagram, by $\sigma_l={\rm diag}(-I_l, I_{n-l},-I_l, I_{n-l})$ for $l=1,\ldots, \left[\frac{n}{2}\right]$; $a_\lambda={\rm diag}(\lambda I_n,\lambda^{-1}I_n)$ and $c_\lambda={\rm diag}(\lambda, I_{n-1},\lambda^{-1},I_{n-1})$ for $\lambda^2\neq0,1$. 
The unipotent ones are those associated with the partitions $(2^{2m},1^{2n-4m})$ for $m=1,\ldots,\left[\frac{n}{2}\right]$ and $(3, 2^{2m},1^{2n-3-4m})$ for $m=1,\ldots,\left[\frac{n}{2}\right]-1$.
\end{theorem}
\pf We view $SO_{2n}(k)$ as the subgroup of $SL_{2n}(k)$ of matrices preserving the bilinear form associated with the matrix $\begin{pmatrix}0&I\cr
I&0\cr\end{pmatrix}$ with respect to the canonical basis of $k^{2n}$. 
We choose $B$ as the subgroup of $G$ of matrices of the form $\begin{pmatrix}A&AX\cr
0&^t\!A^{-1}\cr\end{pmatrix}$ where $A$ is an invertible upper triangular matrix, $^tA^{-1}$ is its inverse transpose and $X$ is a skew symmetric matrix. We fix $T\subset B$ as its subgoup of diagonal matrices.

Let $s\in T$. If $s$ has at least $4$ eigenvalues, up to a Dynkin diagram automorphism $s$ is conjugate to some $r={\rm diag}(\lambda,\mu,\nu,t,\lambda^{-1},\mu^{-1},\nu^{-1},t^{-1})$ for a suitable invertible diagonal submatrix $t$ and some scalars $\lambda$, $\mu$, $\nu$, 
with $\nu$ possibly equal to $\lambda^{-1}$.
Such a matrix is regular in $\langle T,\,X_{\pm\alpha_1},\, X_{\pm\alpha_2}\rangle$ so $\O_s$ is not spherical by Remark \ref{regular}. 

Let us now assume that $s$ has exactly $3$ eigenvalues, say, $\lambda, \lambda^{-1}, 1$ and that the multiplicity of $\lambda$ is greater than $1$. Then $\O_s$ contains a matrix of the
form $r$ as above with $\nu=\lambda^{-1}$ and $\mu=1$ and the above argument shows that $\O_s$ cannot be spherical.
On the other hand, if $s\in T$ has exactly $3$ eigenvalues $\lambda^{\pm1},1$ with the multiplicity of $\lambda$  equal to $1$, up to an automorphism of the Dynkin diagram, $s$ is conjugate to a matrix of the form $c_\lambda={\rm
  diag}(\lambda, I_{n-1},\lambda^{-1}, I_{n-1})$. Its centralizer is equal to the
identity component of the centralizer of the involution $\sigma_1={\rm
  diag}(-1,I_{n-1},-1, I_{n-1})$. Since $G/H$ is spherical if and
only if $G/H^\circ$ is so, $\O_{c_\lambda}$ is spherical. 

Let now $s$ have exactly two eigenvalues. If they are $\pm1$
then $s^2=1$ and $\O_s$ is symmetric. If the eigenvalues are
$\lambda^{\pm1}$ with $\lambda^2\neq1$ we may assume that $s$ is conjugate to $a_\lambda={\rm diag}(\lambda I_n,\lambda^{-1}I_n)$ whose centralizer is independent of $\lambda$ in the given range. Since for $\zeta$ a primitive fourth root of $1$ we
have $a_{\zeta}^2\in Z(G)$, by Remark \ref{same} all those classes are
symmetric spaces, hence spherical. 

Let us now consider unipotent classes in $G$. There might be two classes in $G$ associated by Jordan theory with 
the same partition.
However, two such classes are mapped onto each other by a group automorphism araising from the automorphism of the Dynkin diagram. Thus, one is spherical if and only if the other is so and we will not need to distinguish them.    

Let $u$ be a unipotent element in $G$.  It is well-known that the even terms in its associated partition $\underline{\lambda}$ occur pairwise. Let $\underline{\mu}=(\mu_1,\ldots,\mu_l)$ be obtained taking a representative of each term occurring pairwise in $\underline{\lambda}$ and let $\underline{\nu}=(\nu_1,\ldots,\nu_m)$ be obtained by taking the remaining distinct odd terms so that $2n=2|\underline{\mu}|+|\underline{\nu}|$. A representative $u'$ of $\O_u$ can be taken in the subgroup isomorphic to $SO_{\nu_1+\nu_2}(k)\times\cdots \times SO_{\nu_{m-1}+\nu_m}(k)\times SO_{2\mu_1}(k)\times\cdots\times SO_{2\mu_l}(k)$ obtained by repeatedly immerging  $SO_{2d_1}(k)\times SO_{2d_2}(k)$ into $SO_{2(d_1+d_2)}(k)$ by $\left(\begin{pmatrix}A_1&B_1\\C_1& D_1\end{pmatrix}, \begin{pmatrix}A_2& B_2\\
C_2& D_2\end{pmatrix}\right)\mapsto\begin{pmatrix}
A_1&&B_1\\
&A_2&&B_2\\
C_1&&D_1\\
&C_2&&D_2
\end{pmatrix}$.  The component of $u'$ in $SO_{\nu_i+\nu_{i+1}}(k)$ is
associated with the partition $(\nu_i,\nu_{i+1})$. The component of
$u'$ in $SO_{2\mu_i}(k)$ can be chosen to lie and be regular in the
subgroup isomorphic to $SL_{\mu_i}(k)$ obtained through the immersion
$M\mapsto {\rm diag}(M, ^t\!M^{-1})$. As in type $C_n$ this forces
$\mu_i\leq2$ for all $i$. We shall now show that $\nu_1\leq 3$. Let $\nu_1=2l+1$ and $\nu_{2}=2j-1$ with $l\geq j$. For the component $x$ of $u'$ in $SO_{\nu_1+\nu_2}(k)$ we take
$x=\begin{pmatrix}A\\ AX& ^tA^{-1}\end{pmatrix}$ where $A=\begin{pmatrix}1\\
1&\ddots\\
&\ddots&\ddots\\
&&1&1\\
\end{pmatrix}$ and $AX=\begin{pmatrix}
0_{l-1}\\
&1&1\\
&-1&0\\
&&&&0_{j-1}
\end{pmatrix}$. 
Let us put $\gamma_{l}=\alpha_l+2(\alpha_{l+1}+\cdots+\alpha_{n-2})+\alpha_{n-1}+\alpha_n$. Then the image of $x$ in $G$ 
%
lies in $X_{-\alpha_{1}}\cdots X_{-\alpha_{l+j-1}}X_{-\gamma_l}B$ and it lies in the union of the cells corresponding to
$s_{1}\cdots s_{l-1}s_{l+1}s_{l+2}\cdots s_{l+j-1}s_{\gamma_l}$, $s_{1}\cdots s_{l} s_{l+2}s_{l+3}\cdots s_{l+j-1}s_{\gamma_l}$ and $s_{1}\cdots s_{l+j-1}s_{\gamma_l}$. Since none of the involved elements of $W$ is an involution unless $l+j\leq 2$ we conclude that $\nu_1\leq3$.

Conversely, let $\O_u$ be a unipotent class corresponding to $(2^{2m},1^{2n-4m})$ or to
%
%
$(3, 2^{2m}, 1^{2n-3-4m})$. The matrices in \cite[Theorem 12]{ccc} represent these classes also when ${\rm char}(k)$ is odd and their $B$-orbits satisfy the condition in Theorem \ref{dime}. 

We will show that there is no spherical element with Jordan decomposition $g=su$ with $s\not\in Z(G)$ and $u\neq1$. We might assume that $s=c_\lambda$ or $\sigma_l={\rm diag}(-I_l, I_{n-l},-I_l,I_{n-l})$ because $\dim B=\dim\O_{a_\lambda}<\dim\O_{a_\lambda u}$.  

The subgroup $G_{c_\lambda}^\circ$ is of type $D_{n-1}\times k^*$ and it is generated by $T$ and the root subgroups $X_{\pm \alpha_2},\ldots, X_{\pm\alpha_{n}}$. 
If $u\in G_{c_\lambda}^\circ$ 
corresponds to a partition $\underline{\mu}$ of $2n-2$ with repeated terms, then $\mu=(2,2,\underline{\mu}')$ for some partition $\underline{\mu}'$ of $2n-4$ and $c_\lambda u$ would be conjugate to a matrix of the form 
$\begin{pmatrix}
\lambda\\
1&1\\
&1&1\\
&&&A_1&&&&A_1X_1\\
&&&&\lambda^{-1}&-\lambda^{-1}\\
&&&&&1&-1\\
&&&&&&1\\
&&&&&&&^t\!A_1^{-1}
\end{pmatrix}$
where $\begin{pmatrix}
A_1&A_1X_1\\
&^t\!A_1^{-1}
\end{pmatrix}$ represents a unipotent element associated with $\underline{\mu}'$ and $A_1$ is upper triangular. Such a representative lies in $Bs_1s_2B$. Thus, $u$ has no repeated terms and it must be associated with  $(3,1^{2n-5})$. Then $su$ is conjugate to some 
$x=\begin{pmatrix}
A\\
AX&^tA^{-1}
\end{pmatrix}$ where $A=\begin{pmatrix}
\lambda\\
1&1\\
&1&1\\
&&&I_{n-3}
\end{pmatrix}$ and $AX=\begin{pmatrix}
0\\
&1&1\\
&-1&0\\
&&&0_{n-3}
\end{pmatrix}$ so $x\in X_{-\alpha_1}X_{-\alpha_2}X_{-(\alpha_2+2\alpha_3+\cdots+2\alpha_{n-2}+\alpha_{n-1}+\alpha_n)}$ and it lies in the cell corresponding to $s_1s_2s_{\alpha_2+2\alpha_3+\cdots+2\alpha_{n-2}+\alpha_{n-1}+\alpha_n}$ which is not an involution. Thus, $s=\sigma_l$ for some $l$.
The argument used for $G$ of type $C_n$ can be adapted in order to
show that the partition associated with any component of
$u=(u_1,u_2)\in SO_{2l}(k)\times SO_{2n-2l}(k)=G_{\sigma_l}^\circ$
cannot have repeated terms of size greater than $1$. It follows that
such components can only be associated with $(3,1,\ldots,1)$. Suppose
that $u_2$ is nontrivial so that $n\geq l+2$. We may find a
representative $x$ of $\O_g$
as follows
$x=\begin{pmatrix}
A_1&&&A_1X_1\\
1&1\\
&1&1\\
&&&I_{n-l-2}\\
&&&& ^t\!A_1^{-1}&1\\
&1&1&&&1&-1\\
&-1&0&&&&1\\
&&&&&&&I_{n-l-2}
\end{pmatrix}$ where $\begin{pmatrix}
A_1&A_1X_1\\
&^t\!A_1^{-1}
\end{pmatrix}$ represents $u_1$ and $A_1$ is an upper triangular invertible
matrix. Then  $x$
lies in $TX_{-\alpha_l}X_{-\alpha_{l+1}}X_{-(\alpha_{l+1}+2\alpha_{l+2}+\cdots+2\alpha_{n-2}+\alpha_{n-1}+\alpha_n)}$ and in the cell corresponding to $s_ls_{l+1}s_{\alpha_{l+1}+2\alpha_{l+2}+\cdots+2\alpha_{n-2}+\alpha_{n-1}+\alpha_n}$ leading to a contradiction. Similarly, we may exclude that $u_1$ corresponds to $(3,1,\ldots,1)$, concluding the proof.\hfill$\Box$

\subsection{Type $B_n$}
\begin{theorem}Let $G=SO_{2n+1}(k)$.  The spherical semisimple
  classes in $G$ are represented by $\rho_l={\rm diag}(1,-I_l,
  I_{n-l}, -I_l, I_{n-l})$ for $l=1,\ldots,n$; by $d_\lambda={\rm
  diag}(1,\lambda, I_{n-1},\lambda^{-1},I_{n-1})$ and by $b_\lambda={\rm
  diag}(1,\lambda I_n, \lambda^{-1}I_n)$ with $\lambda^2\neq0,1$.
The unipotent ones are those associated with $(2^{2m},1^{2n+1-4m})$ for $m=1,\ldots,\,[\frac{n}{2}]$ and $(3,2^{2m},1^{2n-2-4m})$ for $m=1,\ldots,[\frac{n-1}{2}]$. The spherical classes that are neither semisimple nor unipotent are represented by $\rho_nu$ where $u\in G_{\rho_n}^\circ\cong SO_{2n}(k)$ is a unipotent element associated with $(2^{2m},1^{2n-4m})$ for $m=1,\ldots,[\frac{n}{2}]$. 
\end{theorem}
\pf  We view $G$ as the subgroup of $SL_{2n+1}(k)$ of matrices preserving the bilinear form associated with $\begin{pmatrix}1&0&0\cr
0&0&I_n\cr
0&I_n&0\cr\end{pmatrix}$ with respect to the canonical basis of $k^{2n+1}$. We fix $B$ to be the subgroup of matrices of the form $\begin{pmatrix}1&0&^t\!\gamma\cr
-A\gamma&A&AX\cr
0&0&^t\!A^{-1}\cr\end{pmatrix}$ where $A$ is an invertible upper triangular matrix, $\gamma$ is a column in $k^n$ and the symmetric part of $X$ is $-\frac{1}{2}\gamma ^t\!\gamma$. We fix $T\in B$ as its subgroup of diagonal matrices. 
 
Let $s\in T$ be a spherical element in $G$.
Adapting the proof of Theorem \ref{dn} by using the
embedding $\iota$ of $SO_{2n}(k)$ into $SO_{2n+1}(k)$ given by
$X\mapsto\begin{pmatrix}1&\\
&X\end{pmatrix}$ we see that $s$ has at most $4$ eigenvalues and one of them is $1$. 
Moreover, if they were $4$, we could take $s={\rm
  diag}(1,t,\lambda,-1,t^{-1},\lambda^{-1},-1)$ for some invertible
diagonal submatrix $t$ and some scalar $\lambda$. Then $s$ would be regular in $\langle T,\,X_{\pm\alpha_{n-1}},\,X_{\pm\alpha_n}\rangle$ so by Remark \ref{regular} this case may not occur. 
%
%
Therefore $s$ has $2$ or $3$ eigenvalues. If $s$ has $2$
eigenvalues it is conjugate to the involution $\rho_l={\rm diag}(1,-I_l, I_{n-l},-I_l,
I_{n-l})$ for some $l$.
If $s$ has $3$
eigenvalues they are $1,\lambda,\lambda^{-1}$. Then 
the multiplicity of $\lambda$ and $1$ cannot be both greater than
$1$ as one can see using the embedding $\iota$
and the element $\iota(r')$ adapting the discussion in Theorem \ref{dn}. Thus, $s$ is conjugate either to $b_\lambda={\rm diag}(1,\lambda
I_n,\lambda^{-1}I_n)$  or to $d_\lambda={\rm
  diag}(1,\lambda,I_{n-1},\lambda^{-1}, I_{n-1})$, for
$\lambda^2\neq1,0$. 
A representative of $\O_{b_\lambda}$ satisfying the condition in
Theorem \ref{dime} is to be found in
\cite[Theorem 15]{ccc} and it is well-defined in odd characteristic,
too. Moreover, $G_{d_\lambda}=G_{\rho_1}^\circ$ hence $b_\lambda$ and
$d_\lambda$ are indeed spherical. 

%

Partitions where even terms occur pairwise parametrize unipotent conjugacy classes in $G$.
Let $u$ be a spherical unipotent element in $G$ associated with
the partition $\underline{\lambda}$. Let $\underline{\mu}$ and $\underline{\nu}$ be constructed as in Theorem \ref{dn} with $2n+1=2|\underline{\mu}|+|\underline{\nu}|$. We may find a representative $u'$ of
$\O_u$ in the subgroup isomorphic to $SO_{\nu_1}(k)\times
SO_{\nu_2+\nu_3}(k)\times\cdots\times SO_{\nu_{m-1}+\nu_m}(k)\times
SO_{2\mu_1}(k)\times\cdots\times SO_{2\mu_l}(k)$ obtained using the
embeddings in the proof of Theorem \ref{dn} and the
embedding of
$SO_{2d_1+1}(k)\times SO_{2d_2}(k)$ into $SO_{2(d_1+d_2)+1}(k)$
given by $\left(\begin{pmatrix}
1&\alpha_1&\beta_1\\
\gamma_1&A_1&B_1\\
\delta_1&C_1&D_1
\end{pmatrix},\begin{pmatrix}
A_2&B_2\\
C_2& D_2
\end{pmatrix}\right)\mapsto\begin{pmatrix}
1&\alpha_1&&\beta_1\\
\gamma_1&A_1&&B_1\\
&&A_2&&B_2\\
\delta_1&C_1&&D_1\\
&&C_2&&D_2
\end{pmatrix}$. The component of $u'$ in $SO_{\nu_i}(k)$ is regular therein,
forcing $\nu_1\leq3$. 
Moreover, the argument used in Theorem \ref{dn} shows that $\mu_i\leq2$ for every
$i$.
%
%
Conversely, for a unipotent 
class associated with $(2^m,1^{2n+1-4m})$ or to
$(3,2^{2m},1^{2n-4m-2})$ the representative in \cite[Theorem 12]{ccc}
is well-defined in odd characteristic and the corresponding 
$B$-orbit satisfies the condition in Theorem \ref{dime}.  

Let $g=su$ be the Jordan decomposition of a spherical element in $G$
with $s,\,u\neq1$. The argument in Theorem \ref{dn}
can be adapted to the present situation, using the
 embedding $\iota$. It follows that for any $\iota(t)\in\O_s$ we have $t\in Z(SO_{2n}(k))$.
%
Therefore $g$ is conjugate to $\rho_l\v={\rm diag}(1,-I_{2n})\v$
for some $\v$ spherical unipotent in $G_{\rho_l}^\circ\cong SO_{2n}(k)$. We claim that
the partition $\underline{\lambda}$ of $2n$ associated with $\v$ has no term equal to $3$.
Otherwise $\underline{\lambda}=(3,\underline{\lambda'},1)$ for some partition
$\underline{\lambda}'$ and $g$ is conjugate to a matrix of the form $$y=\begin{pmatrix}
1&0&1\\
0&-1&0\\
0&1&-1\\
&&&-A&&&-AX\\
1&0&\frac{1}{2}&&-1&-1\\
1&0&\frac{1}{2}&&0&-1\\ 
&&&&&&-^tA^{-1}
\end{pmatrix}
$$
where $\begin{pmatrix}A&AX\\
&^tA^{-1}
\end{pmatrix}$ is a unipotent matrix in $SO_{2n-4}(k)$ corresponding
to $\underline{\lambda}'$ and $A$ is upper triangular. The element $y$ lies in $TX_{-\alpha_1}X_{-(\alpha_2+\cdots+\alpha_n)}B$ and it is contained in the Bruhat cell corresponding to $s_1s_{\alpha_2+\cdots+\alpha_n}$ which is not an involution. Conversely, let $g=\rho_nu$ with $u$ corresponding to $(2^{2m},1^{2n-4m})$. The representative of its class provided in \cite[Theorem 21]{ccc} is well defined in odd characteristic and it allows application of Theorem \ref{dime}.\hfill$\Box$

\subsection{Type $E_6$}
\begin{theorem}\label{e6}Let $G$ be simply connected of type $E_6$. Then the
  spherical classes are either semisimple or unipotent. The
  semisimple ones are symmetric spaces, and, up
  to a central factor, they are represented by $p_1= h_1(-1)h_4(-1)h_6(-1)$ and  $p_{2,c}=h_1(c^2)h_2(c^3)h_3(c^4)h_4(c^6)h_5(c^5)h_6(c^4)$ for
  $c\in k$ with $c^3\neq1,0$. The unipotent ones are those of type $A_1,\,2A_1$ and $3A_1$.
\end{theorem}
\pf Let $s\in T$ be spherical. By
\cite[Theorem 2.15]{Hu-cc} we may choose $s$ so that $G_s$  is
generated by $T$ and $X_{\pm\alpha}$ for $\alpha$ in a subsystem
$\Phi(\Pi)\subset \Phi$ with basis a subset $\Pi$ of $\Delta\cup\{-\beta_1\}$.  
By Theorem \ref{dime}  we have $\dim\O_s\leq \ell(w_0)+{\rm
  rk}(1-w_0)$ and a dimension counting shows that $\Pi$ is one of
the following subsets:
$\Pi_1=\{\alpha_1,\alpha_3,\alpha_4,\alpha_5,\alpha_6,-\beta_1\}$, 
$\Pi_2=\{\alpha_1,\alpha_2,\alpha_4,\alpha_5,\alpha_6,-\beta_1\}$ and 
$\Pi_3=\{\alpha_1,\alpha_2,\alpha_3,\alpha_4,\alpha_6,-\beta_1\}$
of type $A_5\times A_1$; 
$\Pi_4=\{\alpha_1,\alpha_2,\alpha_3,\alpha_4,\alpha_5\}$,
$\Pi_5=\{\alpha_2,\alpha_3,\alpha_4,\alpha_5,-\beta_1\}$, and 
$\Pi_6=\{\alpha_2,\alpha_3,\alpha_4,\alpha_5,\alpha_6\}$
of type $D_5\times k^*$.

Let us put $H_i=\langle T, X_{\pm\alpha},\alpha\in\Pi_i\rangle$ for $i=1,\ldots,6$. The sets $\Pi_i$ for
  $i=1,2,3$ are ${\mathbb R}$-bases for the span of $\Delta$ and
  one may find automorphisms of $\Phi$ mapping $\Pi_i$ for $i=2,3$ to
  $\pm\Pi_1$. 
On the other hand, ${\rm Aut}(\Phi)=\{-w_0\}\ltimes W$ so any element
$s$ whose centralizer is $H_2$ or $H_3$ is $N(T)$-conjugate to an
element whose centralizer is $H_1$. 
The elements $s$ for which $G_s=H_1$ are  $p_1= h_1(-1)h_4(-1)h_6(-1)$ and $zp_1$ for any $z\in Z(G)$. Conjugation by these
  elements is an involution so $\O_{p_1}$ is a symmetric space. 
This completes the analysis for $\Pi_i$ with $i\leq3$.  

The subgroups $H_4$ and
$H_6$ are
$\w_0$-conjugate, so any element whose centralizer is $H_4$ is
$N(T)$-conjugate to an element whose centralizer is $H_6$. Besides,
the automorphism of $\Phi$ defined by
$\alpha_1\mapsto-\beta_1$, $\alpha_2\mapsto \alpha_3$,
$\alpha_3\mapsto\alpha_2$, $\alpha_j\mapsto\alpha_j$ for $j=4,5,6$
maps $\Pi_4$ onto $\Pi_5$ so $H_5$ is $N(T)$-conjugate to $H_4$ and  any element whose centralizer is $H_5$ is
$N(T)$-conjugate to an element whose centralizer is $H_4$. 
The elements whose centralizer is $H_4$ are
$p_{2,c}=h_1(c^2)h_2(c^3)h_3(c^4)h_4(c^6)h_5(c^5)h_6(c^4)$ for $c\in k$ with
  $c^3\neq1,0$. Multiplying $c$ by a third root of unity yields
  the same element multiplied by a central one. Since $p_{2,-1}$ is an involution, $\O_{p_{2,c}}$ is a symmetric space by
  Remark \ref{same}. We claim that $p_{2,c}$ is not conjugate to
  $p_{2,d}$ for $c\neq d$.   
If they were  $G$-conjugate they would be $N(T)$-conjugate by \cite[\S 3.1]{131} so there would
  be $\sigma\in W$ such that $\dot\sigma p_{2,c}\dot\sigma^{-1}=p_{2,d}$.
  Thus, $\sigma$ would stabilize $\Phi(\Pi_4)$ and the restriction of
  $\sigma$ to $\Phi(\Pi_4)$ is an automorphism. It is therefore
  of the form $\tau w$ where $\tau$ acts an automorphism of the Dynkin
  diagram of type $D_5$ and $w$ lies in the Weyl group $W'$ of $H_4$,
  which is contained in $W$. Then $\sigma w^{-1}$ acts on $\Pi_4$ as
  $\tau$. Besides, two
  automorphisms of $\Phi$ coinciding on $\Pi_4$ are equal. Indeed,  
  $\alpha_6$ is the only root $\alpha$ which is orthogonal to $\alpha_j$ for
  $j=1,2,3,4$ and for which $\langle\alpha,\alpha_5\rangle=-1$. It
  follows that  $\sigma w^{-1}$ is either the identity, when
  $\tau=1$,  
 or it is the automorphism mapping
  $\alpha_j$ to $\alpha_j$ for $j=1,3,4$, interchanges $\alpha_2$ and
  $\alpha_5$ and maps $\alpha_6$ to $-\beta_1$. However, one may
  verify that the second possibility cannot happen because such an
  automorphism is equal to
$s_1s_3s_4s_5s_2s_4s_6s_5s_3s_4s_1s_3s_2s_4s_5s_6(-w_0)$, hence it
  does not lie in $W$. Therefore $\tau=1$ and $\sigma=w\in W'$. It is
  not hard to verify that conjugation by the lift in $N(T)$ of an element in
  $W'$
 does not modify $p_{2,c}$ so these elements represent distinct classes.
%

Let $\O$ be a nontrivial spherical unipotent class. Then $\dim\O\leq
\ell(w_0)+{\rm rk}(1-w_0)$ by Theorem \ref{dime} so $\O$ is of type $A_1$, $2A_1$ or
$3A_1$. Conversely, the arguments in \cite[Theorem 13]{ccc} apply in
good characteristic and show that the listed orbits have a
representative whose $B$-orbit satisfies the conditions of Theorem
\ref{dime}. 

A dimension counting together with Lemma \ref{mixed} shows that no class $\O_{su}$ with $s\not\in Z(G)$ and $u\neq1$ can be spherical.\hfill$\Box$

\subsection{Type $E_7$}
\begin{theorem}Let $G$ be simply connected of type $E_7$. The
  spherical classes in $G$ are either semisimple or unipotent. The
  semisimple ones are symmetric and
  are represented by $q_{1}=h_2(\zeta)h_5(-\zeta)h_6(-1)h_7(\zeta)$ where
$\zeta$ is a fixed primitive fourth root of $1$;
$q_2=h_3(-1)h_5(-1)h_7(-1)$; $zq_1$, $zq_2$ for
  $z\in Z(G)$, and 
$q_{3,a}=h_1(a^2)h_2(a^3)h_3(a^4)h_4(a^6)h_5(a^5)h_6(a^4)h_7(a^3)$ with
$a^2\neq1,0$. The unipotent ones are those of type $A_1$, $2A_1$, $(3A_1)'$, $(3A_1)''$ and $4A_1$.
\end{theorem}
\pf Let $s\in T$ be a spherical element. Proceeding as in Theorem \ref{e6}
we may choose $s$ so that $G_s$ is
generated by $T$ and $X_{\pm\alpha}$ for $\alpha\in\Phi(\Pi)$ with
$\Pi$ one of the following subsets of $\Delta\cup\{-\beta_1\}$:
$\Pi_1=\{\alpha_1,\alpha_3,\alpha_4,\alpha_5,\alpha_6,\alpha_7, -\beta_1\}$
of type  $A_7$;  $\Pi_2=\{\alpha_2,\alpha_3,\alpha_4,\alpha_5,\alpha_6,\alpha_7,
-\beta_1\}$, $\Pi_3=\{\alpha_1,\alpha_2,\alpha_3,\alpha_4,\alpha_5,\alpha_7,
-\beta_1\}$
of type $D_6\times A_1$; 
and
$\Pi_4=\{\alpha_1,\alpha_2,\alpha_3,\alpha_4,\alpha_5,\alpha_6\}$ of type
$E_6\times k^*$. 
Let us put $H_i=\langle T, X_{\pm\alpha}, \alpha\in \Pi_i\rangle$.

There is only one element, up to a central one, whose centralizer is $H_1$ and this is
$q_{1}=h_2(\zeta)h_5(-\zeta)h_6(-1)h_7(\zeta)$ where $\zeta$ is a fixed primitive fourth root of $1$. Since $q_1^2=h_2(-1)h_5(-1)h_7(-1)\in Z(G)$ the
corresponding class is symmetric. The root systems generated by $\Pi_2$ and $\Pi_3$ are mapped onto
each other by automorphisms of $\Phi$ which therefore lie in $W$. Thus,
each element whose centralizer is $H_2$ is $N(T)$-conjugate to one
whose centralizer is $H_3$ and it is enough to look at
$\Pi_2$. The elements whose centralizer is $H_2$ are
$q_2=h_3(-1)h_5(-1)h_7(-1)$ and $zq_2$ for $z\in Z(G)$. The corresponding classes
are symmetric. 
The elements whose centralizer is $H_4$ are:
$q_{3,a}=h_1(a^2)h_2(a^3)h_3(a^4)h_4(a^6)h_5(a^5)h_6(a^4)h_7(a^3)$ for
$a^2\neq1,0$. For $\xi$ a primitive fourth root of unity we have
$q_{3,\xi}^2\in Z(G)$ hence all such classes are symmetric
spaces. Multiplication of $q_{3,a}$ by the nontrivial central element
gives $q_{3,-a}$. We claim that $q_{3,a}$ is never conjugate to $q_{3,b}$ for $a\neq
b$. For, if they are $G$-conjugate, they are $N(T)$-conjugate, and there is a $\sigma\in W$
for which $\dot\sigma q_{3,a}\dot\sigma^{-1}=q_{3,b}$. 
Then $\sigma$ preserves $\Phi(\Pi_4)$ and its restriction to it is an automorphism. As in the proof of Theorem \ref{e6} we see that there
is some $w$ in the Weyl group $W'$ of $H_4$ for which
the restriction to $\Phi(\Pi_4)$ of $\sigma w^{-1}\in W$ is an automorphism
of the Dynkin diagram of type $E_6$. However, there is no automorphism of
$\Phi$ whose restriction to $E_6$ is the nontrivial automorphism
because there is no $\alpha\in\Phi$  which is orthogonal to $\alpha_j$ for
$j=2,3,4,5,6$ and for which
$\langle\alpha,\alpha_1\rangle=-1$. Therefore
$\sigma w^{-1}$ is the identity on $\Phi(\Pi_4)$. By uniqueness of the
extension of an automorphism from $E_6$ to $E_7$ we have $w=\sigma$.
It is not hard to verify that conjugation
by lifts in $N(T)$ of elements in $W'$ preserves $q_{3,a}$.

Let $u\neq1$ be a spherical unipotent element. Then
$\dim\O_u\leq\dim B$ so $\O_u$ is either of type $rA_1$ for some $r$
or of type $A_2$. In the latter case, $u$ would be distinguished in a
standard Levi subgroup of type $A_2$. Since in type $A_n$
distinguished elements are regular, this case cannot occur by Remark \ref{regular}. 
The arguments in \cite[Theorem 13]{ccc} apply also in good
characteristic and they show that for all unipotent classes
of type $rA_1$ there is a representative whose $B$-orbit satisfies the
condition in  Theorem \ref{dime}.

We claim that there is no spherical element with Jordan decomposition
$g=su$ with $s\not\in Z(G)$ and $u\neq1$. Indeed,  $\O_s$ would 
be spherical and $u$ would be spherical in $G_s^\circ$. 
A dimensional argument shows that this could be possible only if $s\in \O_{q_2}$
 and $u$ is nontrivial only in the component of type $A_1$
 in $G_s$. Then $g$ would be conjugate to
 $q_2x_{-\alpha_7}(1)$. 
Conjugation by $x_{-\alpha_6}(1)$ and Chevalley's commutator formula would give 
$q_2 x_{-\alpha_6}(a)x_{-\alpha_7}(1)x_{-\alpha_6-\alpha_7}(b)\in
\O_g$ for some nonzero $a,\,b\in k$. Conjugating by a suitable element
in $X_{-\alpha_6-\alpha_7}$ we could get rid of the term in
$X_{-\alpha_6-\alpha_7}$ obtaining an element in $\O_g\cap Bs_6s_7
B$. By Theorem \ref{invo} the class $\O_g$ cannot be spherical.\hfill$\Box$

\subsection{Type $E_8$}
\begin{theorem}Let $G$ be of type $E_8$. The spherical classes are either semisimple or unipotent. The semisimple ones
 are symmetric and they are represented by $r_1=h_2(-1)h_3(-1)$ and
 $r_2=h_2(-1)h_5(-1)h_7(-1)$. The unipotent ones are those
 of type $A_1$, $2A_1$, $3A_1$ and $4A_1$.
\end{theorem}
\pf Let $s\in T$ be a spherical element. Since
$\dim\O_s\leq\dim B$, up to conjugation in $N(T)$, the centralizer $G_s$ is
generated by $T$ and by the $X_{\pm\alpha}$ for
$\alpha$ in a subsystem generated either by $\{\alpha_2,\alpha_3,\alpha_4,\alpha_5,\alpha_6,\alpha_7,\alpha_8,-\beta_1\}$
of type $D_8$ or by $\{\alpha_1,\alpha_2,\alpha_3,\alpha_4,\alpha_5,\alpha_6,\alpha_7,-\beta_1\}$
of type $E_7\times A_1$. 
Then $s$ is conjugate either to $r_1=h_2(-1)h_3(-1)$ or to $r_2=h_2(-1)h_5(-1)h_7(-1)$.
Since $r_1^2=r_2^2=1$ the corresponding classes are symmetric spaces.\\ 
Let $\O$ be a nontrivial spherical unipotent class. Then
$\dim \O\leq\dim B$ so $\O$ is either of type $rA_1$ for some $r$, or
it is of type $A_2$. The latter case is excluded as we did for $G$ of
type $E_7$. Conversely, the arguments in \cite[Theorem 13]{ccc} apply in
good characteristic and they show that for each orbit of type $rA_1$ in $G$ we may find a representative whose $B$-orbit satisfies the condition in Theorem \ref{dime}. 

%

We claim that there is no spherical element with Jordan decomposition $g=su$ 
with $s,u\neq 1$. Indeed, by dimensional reasons $s$ would be
conjugate to $r_2$ and $u$ would lie in the component of
type $A_1$ in $G_s$.
In other words, we may assume $g=r_2 x_{-\beta_1}(1)$. Let
$\gamma=\beta_1-\alpha_8$. Conjugation of $g$ by $\s_\gamma$ gives $t
x_{-\alpha_8}(a)\in\O_g$ for some nonzero $a\in k$ and some $t\in
T$. Since $r_2$ does not commute with
$X_{\pm(\beta_1-\alpha_8-\alpha_7)}$ and
$s_\gamma(\alpha_7+\alpha_8-\beta_1)=\alpha_7$ the element $t$ does
not commute with $X_{\pm\alpha_7}$. Since
$s_\gamma(\alpha_7+\alpha_8)=\alpha_7+\alpha_8$ and $r_2$ does not
commute with $X_{\pm (\alpha_7+\alpha_8)}$ the same holds for
$t$. Then conjugation of $tx_{-\alpha_8}(a)$ by $x_{-\alpha_7}(1)$
would give $t
x_{-\alpha_7}(b)x_{-\alpha_8}(a)x_{-\alpha_7-\alpha_8}(c)\in\O_g$ for
some nonzero $b,c\in k$. Conjugation by a suitable element in
$X_{-\alpha_7-\alpha_8}$ would yield an element in $\O_g\cap Bs_7s_8
B$ concluding the proof. \hfill$\Box$

\subsection{Type $F_4$}
\begin{theorem}\label{f4}Let $G$ be of type $F_4$.  The spherical semisimple classes are symmetric spaces and they are represented by $f_1=h_{\alpha_2}(-1)h_{\alpha_4}(-1)$ and $f_2=h_{\alpha_3}(-1)$.
The spherical unipotent ones are those of type
$rA_1+s\tilde{A}_1$ for $r,s\in\{0,1\}$. Furthermore, there is a
spherical class that is neither semisimple nor unipotent and
it is represented by $f_2x_{\beta_1}(1)$. 
\end{theorem}
\pf Let $s\in T$ be a spherical semisimple element in $G$. 
A  dimension counting shows that $G_s$ is $N(T)$-conjugate to the
subgroup generated by $T$ and the root subgroups in a subsystem
generated either by 
 $\Pi_1=\{-\beta_1,\,\alpha_2,\alpha_3,\,\alpha_4\}$ or by
$\Pi_2=\{\alpha_1,\alpha_2,\,\alpha_3,\,-\beta_1\}$. They correspond
to $f_1=h_{\alpha_2}(-1)h_{\alpha_4}(-1)$ and $f_2=h_{\alpha_3}(-1)$,
respectively.

Let $\O$ be a nontrivial spherical unipotent class in
$G$. Then $\dim\O\leq\dim B$ so $\O$ is either of type $A_1,
\tilde{A}_1$ or $A_1+\tilde{A}_1$. Conversely, the arguments in \cite[Theorem
  13]{ccc} hold in good characteristic and they show that
Theorem \ref{dime} applies to all classes of type $rA_1+s\tilde{A}_1$ in $G$.
%
%

Let $g=su$ be the Jordan decomposition of a spherical element with $s,\,u\neq1$. Since $\dim\O_{f_1}=\dim B$ we
might assume $s=f_2$. Besides, $G_{f_2}$ is a reductive group of type
$B_4$. A dimensional argument shows that $u$ lies in the minimal unipotent
class in $G_{f_2}$ so we may assume $g=f_2
x_{2\alpha_1+3\alpha_2+4\alpha_3+2\alpha_4}(1)$. We have
$\dim\O_g=\dim B$. The proof in \cite[Theorem 23]{ccc} 
contains an incorrect argument which we rectify here. 

The element $f_2=h_{\alpha_3}(-1)$ lies in the subgroup $G_1=\langle X_{\pm\alpha_i},\;i=2,3,4\rangle$ of type $C_3$. By looking at the centralizer of $f_2$ in $G_1$ we see that,  up to an element in $Z(G_1)$, the conjugacy class of $f_2$ in $G_1$ corresponds to $\sigma_1$ with notation as in Theorem \ref{cn}. By \cite[Theorem 15]{ccc}, the class of $\sigma_1$ in $G_1$
has a representative in $s_{4}s_{\alpha_2+2\alpha_3+\alpha_4}T$ when
$k={\mathbb C}$. The
same matrix represents the class in good characteristic. Besides, $G_1$ centralizes $X_{\pm\beta_1}$, so $f_2X_{-\beta_1}$ can
be represented by an element $z\in
s_{4}s_{\alpha_2+2\alpha_3+\alpha_4}TX_{-\beta_1}\subset
X_{\beta_1}w_0s_2TX_{\beta_1}$. Conjugating $z$ by $\s_2\s_1$ we obtain an element $z'\in B w_0 s_1 B\cap\O_g$. Thus, $w_{\O_g}\geq w_0s_2$ and $w_{\O_g}\geq w_0s_1$, forcing $w_0=w_\O$. Then $\O_g$ has a representative whose $B$-orbit satisfies  the condition in Theorem \ref{dime}.\hfill$\Box$     

\subsection*{Acknowledgements}
The author wishes to thank Mauro Costantini for pointing out an
incorrect argument in the discussion of type $F_4$ in \cite{ccc} and in a previous version of this manuscript.
It is rectified in Theorem \ref{f4}.

\end{document}